\documentclass[11pt, reqno]{amsart}
\usepackage{mathrsfs}
\usepackage{mathrsfs}

\usepackage{color}
\usepackage[unicode]{hyperref}
\hypersetup{
	colorlinks = true,%
	citecolor = [rgb]{0.0,0.0,0.9},
	filecolor=black,%
	linkcolor = [rgb]{0.65,0.0,0.0},%
	anchorcolor = red,
	pagecolor = red,
	urlcolor= [rgb]{0.65,0.0,0.0}
}

\newtheorem{teor}{Theorem}[section]

\newtheorem{remark}{Remark}[section]

\newcommand{\R}{\Bbb{R}}

\newcommand{\n}{\noindent}

\begin{document}

\title[Second-order Sobolev inequalities on Riemannian manifolds]{Second-order Sobolev inequalities on a class of Riemannian manifolds with nonnegative Ricci curvature}


\author{Ezequiel Barbosa}
\address{Departamento de Matem\'{a}tica, Universidade Federal de Minas Gerais, Caixa Postal 702, 30123-970, Belo Horizonte, MG, Brazil.}
\email{ezequiel@mat.ufmg.br}

\author{Alexandru Krist\'aly}
\address{Department of Economics, Babe\c s-Bolyai University, 400591 Cluj-Napoca,
	Romania \& Institute of Applied Mathematics, \'Obuda
	University, 1034 
	Budapest, Hungary.}
\email{alex.kristaly@econ.ubbcluj.ro; alexandrukristaly@yahoo.com}
\thanks{E. Barbosa is supported by CNPq-Brazil. A. Krist\'aly is supported by a grant of the Romanian National Authority for Scientific Research,
CNCS-UEFISCDI, project no. PN-II-ID-PCE-2011-3-0241. 
}




\begin{abstract}
	Let $(M,g)$ be an $n-$dimensional  complete open  Riemannian ma\-ni\-fold with nonnegative Ricci curvature verifying  $\rho\Delta_g \rho \geq n- 5\geq 0$, where $\Delta_g$ is the Laplace-Beltrami operator on $(M,g)$ and $\rho$ is the distance function from a given point. If $(M,g)$ supports a second-order Sobolev inequality with a constant $C>0$ close to the optimal  constant $K_0$ in the second-order Sobolev inequality in $\mathbb R^n$, we show that a  global volume non-collapsing property holds on $(M,g)$. The latter property together with a Perelman-type construction established by  Munn (J. Geom. Anal., 2010)  provide several rigidity results in terms of the higher-order homotopy groups of $(M,g)$. Furthermore, it turns out that $(M,g)$ supports the second-order Sobolev inequality with the constant $C=K_0$ if and only if  $(M,g)$ is isometric to the Euclidean space $\mathbb R^n$. 
\end{abstract}

\maketitle

\vspace{-0.7cm}
\section{Introduction}
It is  well known that the validity of first-order Sobolev inequalities on Riemannian manifolds strongly depend on the curvature; this is a rough conclusion of the famous AB-program initiated by Th. Aubin in the seventies, see the monograph of Hebey \cite{Hebey} for a systematic presentation.  To be more precise, let $(M,g)$ be an $n-$dimensional complete Riemannian manifold, $n\geq 3,$ and consider for some $C>0$ the  first-order Sobolev inequality  
$$\left(\int_M |u|^{2^*}{\rm d}v_g\right)^\frac{2}{2^*}\leq C \int_M |\nabla_g u|^2{\rm d}v_g,\ \forall u\in C^{\infty}_0(M), \eqno{({\bf FSI})_C}$$
where $2^*=\frac{2n}{n-2}$ is the first-order critical Sobolev exponent, and  ${\rm d}v_g$ and $\nabla_g$ denote the canonical volume form and gradient on $(M,g)$, respectively. On one hand, inequality $({\bf FSI})_C$ holds on any $n-$dimensional Cartan-Hadamard manifold $(M,g)$ (i.e., simply connected, complete Riemannian manifold with nonpositive sectional curvature) with the optimal Euclidean constant $C=c_0=\left[\pi n(n-2)\right]^{-1}\left(
{\Gamma(n)}/{\Gamma(\frac{n}{2})} \right)^{{2}/{n}}$ whenever the Cartan-Hadamard conjecture holds on $(M,g)$, e.g., $n\in \{3,4\}$. On the other hand, due to Ledoux \cite{Ledoux1}, if $(M,g)$ has nonnegative Ricci curvature, inequality $({\bf FSI})_{c_0
}$ holds if and only if $(M,g)$ is isometric to the  Euclidean space $\mathbb R^n$. Further  first-order Sobolev-type inequalities on Riemannian/Finsler manifolds can be found in 
 Bakry, Concordet and  Ledoux \cite{Ledoux2},  
Druet,  Hebey and  Vaugon 
 \cite{HebeyDruetVaugon}, do Carmo and Xia \cite{CarmoXia},   Xia \cite{Xia1}-\cite{Xia3},  Krist\'aly \cite{K-JGA}; moreover, similar Sobolev inequalities are also considered  on 'nonnegatively' curved metric measure spaces  formulated in terms of the Lott-Sturm-Villani-type curvature-dimension condition or the Bishop-Gromov-type doubling measure condition, see Krist\'aly \cite{Kristaly-Calculus-2016} and  Krist\'aly and Ohta \cite{KO-Math_Annal}. 

With respect to first-order Sobolev inequalities, much less is know about higher-order Sobolev inequalities on curved spaces. The first  studies concern the AB-program for Paneitz-type operators on compact Riemannian manifolds, see  Djadli, Hebey and Ledoux \cite{DjadliHebeyLedoux}, Hebey \cite{Hebey-JGA} and Biezuner and Montenegro \cite{Biezuner-Montenegro}. Recently, Gursky and Malchiodi \cite{GM} studied strong maximum principles for Paneitz-type operators on complete Riemannian manifolds with semi-positive $Q-$curvature and nonnegative scalar curvature. 

The aim of the present paper is to establish rigidity results on Riemannian manifolds with nonnegative Ricci curvature supporting  second-order Sobolev inequalities.  In order to present our results, let  $(M,g)$ be an $n$-dimensional complete open Riemannian
manifold, $n\geq 5$,  $B(x,r)$
be the geodesic ball with center $x\in M$ and radius $r>0,$ and
${\rm vol}_g[B(x,r)]$ be the volume of $B(x,r)$. We say that $(M,g)$ supports the {\it second-order Sobolev inequality  for  $C>0$} if  
$$
\left(\int_M|u|^{2^{\sharp}}{\rm d}v_g\right)^{\frac{2}{2^{\sharp}}}\leq
C\int_M\left(\Delta_gu\right)^2{\rm d}v_g,\  \ \ \forall u\in C^{\infty}_0(M), \eqno{({\bf SSI})_C}
$$
where $2^{\sharp} = \frac{2n}{n - 4}$ is the second-order critical
Sobolev exponent, and 
$\Delta_g$ is the Laplace-Beltrami
operator on $(M,g)$. Note that the Euclidean space $\mathbb R^n$ supports $({\bf SSI})_{K_0}$ for  
\begin{equation}\label{1a}
K_0=\left[\pi^2n(n-4)(n^2-4)\right]^{-1}\left(
\frac{\Gamma(n)}{\Gamma(\frac{n}{2})} \right)^{{4}/{n}}\,.
\end{equation} Moreover, $K_0$ is optimal, see Edmunds, Fortunato and Janelli \cite{Fortunato}, Lieb \cite{Lieb} and
Lions \cite{Lions}, and the unique class of extremal functions is 
$$u_{\lambda,x_0}(x)=\left(\lambda+|x-x_0|^2\right)^{\frac{4-n}{2}},\ x\in \mathbb R^n,$$
where $\lambda>0$ and $x_0\in \mathbb R^n$ are arbitrarily fixed.

 To state our results, we need a technical assumption on the manifold $(M,g)$; namely, if  $\rho$ is the distance function on
$M$ from a given point $x_0\in M$, we say that $(M,g)$ satisfies the {\it distance Laplacian
growth condition} if
\[
\rho\Delta_g \rho\geq n-5\,.
\]

Now, our main result reads as follows. 
\begin{teor}\label{main}
Let $n\geq5$ and $(M, g)$  be an $n$-dimensional complete open
Riemannian manifold with nonnegative Ricci curvature which satisfies 
the distance Laplacian growth condition. Assume that $(M,g)$ supports the second-order Sobolev inequality $({\bf SSI})_C$ for some $C>0$. Then the following properties hold: 

\begin{itemize}
	\item[(i)] $C\geq K_0;$
	\item[(ii)] if in addition $C\leq \frac{n+2}{n-2}
	K_0$, then we have the  global volume non-collapsing property \[
	{\rm vol}_g[B(x,r)]\geq (C^{-1}K_0)^{\frac{n}{4}}\omega_n r^n \quad for\ all\
	r>0, x\in M\,,
	\]
	where $\omega_n$ is the volume of the $n-$dimensional Euclidean unit ball. 
\end{itemize}
\end{teor}

\begin{remark}\rm 
The distance
	Laplacian growth condition on $(M,g)$ is indispensable in our argument which shows the genuine second-order character of the studied problem. We notice that the counterpart of this condition in the first-order Sobolev inequality $({\bf FSI})_C$ is the validity of an eikonal inequality $|\nabla_g \rho|\leq 1$ a.e. on $M$, which trivially holds on any complete Riemannian manifold (and any metric measure space with a suitable derivative notion).  Further comments on this condition will be given in Section \ref{section-final}.
\end{remark}

Having the global volume non-collapsing  property of geodesic balls of $(M,g)$ in Theorem \ref{main} (ii), we shall prove that once $ C>0$ in $({\bf SSI})_C$ is closer and closer to
the optimal Euclidean constant $K_0$, the Riemannian manifold
$(M,g)$ approaches topologically more and more to the Euclidean
space $\mathbb R^n.$ 
To describe quantitatively this phenomenon, we recall the construction of Munn \cite{Munn-JGA} 
based on the double induction argument of Perelman
\cite{Perelman}. In fact, Munn determined explicit lower bounds for the volume growth of the geodesic balls in terms of certain constants which
guarantee the triviality of the $k$-th
homotopy group $\pi_k(M)$ of $(M,g).$ 
More precisely, let $n\geq 5$ and for $k\in\{1,...,n\}$, let us denote by $\delta_{k,n}>0$  the
smallest positive solution to the equation
$$10^{k+2}C_{k,n}(k)s\left(1+\frac{s}{2k}\right)^k=1$$ in the variable
$s>0$,  where
\begin{equation*}
C_{k,n}(i)= \left\{
\begin{array}{lll}
1 & \mbox{if} & i=0, \\
3+10C_{k,n}(i-1)+(16k)^{n-1}(1+10C_{k,n}(i-1))^n & \mbox{if} & i\in
\{1,...,k\}.
\end{array}%
\right.
\end{equation*} We now consider
the smooth, bijective and increasing function
$h_{k,n}:(0,\delta_{k,n})\to (1,\infty)$ defined by
$$h_{k,n}(s)=\left[1-10^{k+2}C_{k,n}(k)s\left(1+\frac{s}{2k}\right)^k\right]^{-1}.$$
For every $k\in \{1,...,n\},$ let
\begin{equation*}
\alpha_{MP}(k,n)= \left\{
\begin{array}{lll}
1-\left[1+\frac{2}{h_{1,n}^{-1}(2)}\right]^{-1} & \mbox{if} & k=1, \\
1-\left[1+\left( \frac{1+...+\frac{h_{k-1,n}^{-1}(1+\frac{\delta_{k,n}}{2k})}{2(k-1)}}{h_{1,n}^{-1}\left(1+...+\frac{h_{k-1,n}^{-1}(1+\frac{\delta_{k,n}}{2k})}{2(k-1)}\right)}\right)^n\right]^{-1} & \mbox{if} &
k\in \{2,...,n\},
\end{array}%
\right.
\end{equation*}
 be the so-called {\it Munn-Perelman
	constant},  see Munn \cite[Tables 4 and 5, p. 749-750]{Munn-JGA}.  

%



Following the idea from Krist\'aly \cite{Kristaly-Calculus-2016}, our quantitative result reads as follows: 
\begin{teor}\label{theorem-quantitative}
	Under the same assumptions as in Theorem \ref{main}, we have 
	\begin{itemize}
		\item[{\rm (i)}] if  $C\leq \frac{n+2}{n-2}
		K_0$, the order of the fundamental group $\pi_1(M)$ is bounded
		above by $\left(
		\frac{ C}{K_0} \right)^\frac{n}{4}$ $($in particular, if $C<2^\frac{4}{n}K_0$, then $M$ is simply connected$);$
		\item[{\rm (ii)}] if $ C < \alpha_{MP}(k_0,n)^{-\frac{4}{n}}K_0$ for some $k_0\in \{1,...,n\}$ then
		$\pi_1(M)=...=\pi_{k_0}(M)=0;$
		\item[{\rm (iii)}] if $ C < \alpha_{MP}(n,n)^{-\frac{4}{n}}K_0$  then
		$M$ is contractible$;$
		\item[{\rm (iv)}]
		$ C=K_0$   if and only if $(M,g)$ is isometric to the Euclidean space $\mathbb
		R^n.$
	\end{itemize}
\end{teor}

\section{Proof of Theorems \ref{main}\&\ref{theorem-quantitative}}

Throughout this section, we assume the hypotheses of Theorem \ref{main} are verified, i.e.,  $(M, g)$  is an $n$-dimensional complete open
Riemannian manifold with nonnegative Ricci curvature which satisfies 
the distance Laplacian growth condition and supports the second-order Sobolev inequality $({\bf SSI})_C$ for $C>0$.

(i) The inequality $C\geq K_0$ follows in a similar way as in Djadli, Hebey and Ledoux \cite[Lemmas 1.1\&1.2]{DjadliHebeyLedoux} by using a geodesic, normal coordinate system at a given point $x_0\in M$. 

(ii) 
Before starting the proof explicitly,  we notice that one can assume that $C>K_0$; otherwise, if $C=K_0$ then we can assume that $({\bf SSI})_C$ holds with $C=K_0+\varepsilon$, where $\varepsilon>0$ is arbitrarily small, and then letting $\varepsilon\to 0$. Now, we split the proof into five steps. 

{\bf Step 1.} {\it ODE via the Euclidean optimizer.} 
 We consider the function $G:(0,\infty)\rightarrow \R$ defined by
 \[
 G(\lambda)=\int_{\R^n}\frac{{\rm d}x}{\left( \lambda+|x|^2
 	\right)^{n-2}}\,.
 \]
The layer cake representation  shows that for every $\lambda>0$, 
\begin{equation}\label{G-kifejtve}
G(\lambda)=2(n-2)\omega_n\int
_0^{\infty}\frac{t^{n+1}}{\left( \lambda+t^2 \right)^{n-1}}{\rm d} t=\frac{2^{4-n}\pi^\frac{n+1}{2}}{(n-4)\Gamma(\frac{n-1}{2})}\lambda^\frac{4-n}{2}.
\end{equation}
Clearly, $G$ is smooth on $(0,\infty)$. 

We recall now by (\ref{1a}) that 
$$
\left(\int_{\R^n}|u_\lambda|^{2^{\sharp}}{\rm d}x\right)^{\frac{2}{2^{\sharp}}}=
K_0\int_{\R^n}\left(\Delta u_\lambda\right)^2{\rm d}x\,,
$$
where  $$u_\lambda(x)=(\lambda+|x|^2)^\frac{4-n}{2},\ x\in \mathbb R^n,$$
and $\lambda>0$ is arbitrarily fixed.
%
%
%
%
%
%
%
 In terms of the function $G$, the above equality can be rewritten as 
$$
\left(\frac{G''(\lambda)}{(n-2)(n-1)} \right)^{\frac{n-4}{n}} =
K_0(n-4)^2\left\{4G(\lambda)-4\lambda G'(\lambda)+\frac{n-2}{n-1}\lambda^2 G''(\lambda)\right\}\,.
$$
By introducing the function $$G_0(\lambda)=\left(\frac{K_0}{C}\right)^\frac{n}{4}( G(\lambda)-\lambda G'(\lambda)),\ \lambda>0,$$  the latter relation is equivalent to  the ODE
\begin{equation}\label{6}
\left(-\frac{ G_0'(\lambda)}{\lambda(n-2)(n-1)} \right)^{\frac{n-4}{n}} =
C(n-4)^2\left\{4 G_0(\lambda)-\frac{n-2}{n-1}\lambda  G_0'(\lambda)\right\}\,, \ \ \lambda>0.
\end{equation}

%
%
%
%
%
%
%
%
%
%
%
%
%

{\bf Step 2.} {\it ODI via  $({\bf SSI})_C$.} 
Let $x_0\in M$ be the point for which the distance Laplacian growth condition holds and let $F:(0,\infty)\rightarrow \R$ be defined by
\[
F(\lambda)=\int_M\frac{{\rm d}v_g}{\left( \lambda+\rho^2
\right)^{n-2}}\,.
\]
Since $(M,g)$ has nonnegative Ricci curvature, the Bishop-Gromov comparison theorem asserts that ${\rm vol}_g[B(x_0,t)]\leq \omega_nt^n$ for every $t>0$; thus, by the layer cake representation and a change of variables, it turns out that 
\begin{eqnarray}\label{F-null-rendu}
F(\lambda)&=& \int _0^{\infty}{\rm vol}_g\left \{x\in M:
\frac{1}{\left( \lambda+\rho(x)^2 \right)^{n-2}}>s \right \}{\rm d}s\,\nonumber\\
&=&{2}{(n-2)}\int _0^{\infty}{\rm vol}_g[B(x_0,t)]\frac{t}{\left( \lambda+t^2 \right)^{n-1}}{\rm d}t\\
&\leq&{2}{(n-2)}\omega_n\int _0^{\infty}\frac{t^{n+1}}{\left(
\lambda+t^2 \right)^{n-1}}{\rm d}t\,\nonumber\\
&=&G(\lambda). \nonumber
\end{eqnarray}

\n Thus $0< F(\lambda)<\infty$ for every $\lambda>0$, and $F$ is
smooth. In a similar way, 
\begin{eqnarray}\label{F-elso-rendu}
F'(\lambda)=-{2}{(n-2)(n-1)}\int _0^{\infty}{\rm vol}_g[B(x_0,t)]\frac{t}{\left( \lambda+t^2 \right)^{n}}{\rm d}t,
\end{eqnarray}
$$F''(\lambda)={2}{(n-2)(n-1)n}\int _0^{\infty}{\rm vol}_g[B(x_0,t)]\frac{t}{\left( \lambda+t^2 \right)^{n-1}}{\rm d}t,$$
%
%
and for every $\lambda>0$,
\begin{equation}\label{masodrendu}
-\infty<G'(\lambda)\leq F'(\lambda)<0\ \ {\rm and}\ \ 0< F''(\lambda)\leq G''(\lambda)<\infty.
\end{equation}

Let $\lambda>0$ be fixed; we observe that  the function $$w_{\lambda}=
\left(\lambda+\rho^2 \right)^{\frac{4-n}{2}}$$ can be approximated by elements from $C_0^\infty(M)$; in particular, by using an approximation procedure, one can use the function $w_{\lambda}$ as a test-function in $({\bf SSI})_C$. Accordingly, 
\begin{equation}\label{laplace-inequality}
\left(\int_M|w_\lambda|^{2^{\sharp}}{\rm d}v_g\right)^{\frac{2}{2^{\sharp}}}\leq
C\int_M\left(\Delta_gw_\lambda\right)^2{\rm d}v_g,\ \forall \lambda>0.
\end{equation}
A chain rule and the eikonal equation $|\nabla_g \rho|=1$ shows that
$$
\left(\Delta_gw_{\lambda}\right)^2=(n-4)^2(\lambda+\rho^2)^{-n}\left(
\lambda + (3-n)\rho^2 + (\lambda+\rho^2)\rho\Delta_g\rho
\right)^2\,.
$$

 Since the Ricci curvature is nonnegative on $(M,g)$, we first have the distance Laplacian comparison  $\rho\Delta_g\rho\leq n-1$. Thus, 
\begin{equation}\label{mainE1}
\lambda + (3-n)\rho^2 + (\lambda+\rho^2)\rho\Delta_g\rho\leq
2\rho^2+n\lambda\,,\ \forall \lambda>0.
\end{equation}
On the other hand, by the distance Laplacian growth condition, i.e.,  $\rho\Delta_g\rho\geq n-5$, we obtain that
\begin{equation}\label{mainE2}
-(2\rho^2+n\lambda)\leq \lambda + (3-n)\rho^2 +
(\lambda+\rho^2)\rho\Delta_g\rho\,,\ \forall \lambda>0.
\end{equation} 
 Consequently, by (\ref{mainE1})
and (\ref{mainE2}), we have that
$$
|\lambda + (3-n)\rho^2 + (\lambda+\rho^2)\rho\Delta_g\rho|\leq
2\rho^2+n\lambda\,,\ \forall \lambda>0.
$$

\n Thus, it turns out that 
\[
\left(\Delta_gw_{\lambda}\right)^2\leq(n-4)^2(\lambda+\rho^2)^{-n}\left(2\rho^2+n\lambda\right)^2\,.
\]

According to the latter estimate, relation (\ref{laplace-inequality}) can be written in terms of the function $F$ as 
$$
\left(\frac{F''(\lambda)}{(n-2)(n-1)} \right)^{\frac{n-4}{n}} \leq
C(n-4)^2\left\{4F(\lambda)-4\lambda F'(\lambda)+\frac{n-2}{n-1}\lambda^2 F''(\lambda)\right\}\,.
$$
By defining the function $$F_0(\lambda)= F(\lambda)-\lambda F'(\lambda),$$  the latter relation is equivalent to the ordinary differential inequality \begin{equation}\label{633}
\left(-\frac{ F_0'(\lambda)}{\lambda(n-2)(n-1)} \right)^{\frac{n-4}{n}} \leq
C(n-4)^2\left\{4 F_0(\lambda)-\frac{n-2}{n-1}\lambda  F_0'(\lambda)\right\}\,, \ \ \lambda>0.
\end{equation}

{\bf Step 3.} {\it Comparison of $G$ and $F$ near the origin.} We claim that
\[
\liminf\limits_{\lambda\rightarrow 0}\frac{F(\lambda)-\lambda
	F'(\lambda)}{G(\lambda)-\lambda G'(\lambda)}\geq 1\,.
\]
To see this, fix $\varepsilon >0$ arbitrarily small. Since
 \[
 \lim \limits_{t\rightarrow 0}\frac{{\rm vol}_g[B(x_0,t)]}{\omega_n t^n}=1\,,
 \]
 \n there exists a $\delta>0$ such that ${\rm vol}_g[B(x_0,t)]\geq
 (1-\varepsilon)\omega_n t^n$ for all $t\in(0, \delta]$. Thus, by (\ref{F-null-rendu}) and (\ref{F-elso-rendu}), we have
 \begin{eqnarray*}
 	F(\lambda)&\geq& {2}{(n-2)}\int _0^{\delta}{\rm vol}_g[B(x_0,t)]\frac{t}{\left( \lambda+t^2 \right)^{n-1}}{\rm d}t\\
 	&\geq& {2}{(n-2)}\omega_n(1-\varepsilon)\int_{0}^{\delta}\frac{t^{n+1}}{(\lambda+t^2)^{n-1}}{\rm d}t\\
 	&=&{2}{(n-2)}\omega_n\lambda^\frac{4-n}{2}(1-\varepsilon)\int_{0}^{\delta\lambda^{-\frac{1}{2}}}\frac{s^{n+1}}{(1+s^2)^{n-1}}{\rm d}s,
 \end{eqnarray*}
 \n and
 \begin{eqnarray*}
 	-\lambda F'(\lambda)&\geq& {2}{(n-2)(n-1)}\lambda\int _0^{\delta}{\rm vol}_g[B(x_0,t)]\frac{t}{\left( \lambda+t^2 \right)^{n}}{\rm d}t\\
 	&\geq & {2}{(n-2)(n-1)}\omega_n\lambda(1-\varepsilon)\int _0^{\delta}\frac{t^{n+1}}{\left( \lambda+t^2 \right)^{n}}{\rm d}t\\
 	&=&{2}{(n-2)(n-1)}\omega_n\lambda^\frac{4-n}{2}(1-\varepsilon)\int_0^{\delta\lambda^{-\frac{1}{2}}}\frac{s^{n+1}}{(1+s^2)^{n}}{\rm d}s.
 \end{eqnarray*}
 
  Combining this estimates with  relation (\ref{G-kifejtve}), we obtain 
 \[
 \liminf\limits_{\lambda\rightarrow 0}\frac{F(\lambda)-\lambda
 	F'(\lambda)}{G(\lambda)-\lambda G'(\lambda)}\geq
 1-\varepsilon\,.
 \]
 
 \n Letting $\varepsilon\rightarrow 0 $, we get the required claim.

{\bf Step 4.} {\it Global comparison of $G_0$ and $F_0$.} We claim that 
\begin{equation}\label{F0G0-global}
F_0(\lambda)\geq G_0(\lambda),\ \forall \lambda>0.
\end{equation}
 First of all, by Step 3 and the fact that $C>K_0$, we have 
\begin{eqnarray*}
	\liminf\limits_{\lambda\rightarrow 0}\frac{F_0(\lambda)}{G_0(\lambda)}&=&
	\left(\frac{C}{K_0}\right)^{\frac{n}{4}}\liminf\limits_{\lambda\rightarrow 0}\frac{F(\lambda)-\lambda F'(\lambda)}{G(\lambda)-\lambda G'(\lambda)}\\
	&\geq& \left(\frac{C}{K_0}\right)^{\frac{n}{4}}\\&>& 1\,.
\end{eqnarray*}
%
Thus, for sufficiently small $\delta_0>0,$ one has
\begin{equation}\label{F0G0-global-1}
F_0(\lambda)\geq G_0(\lambda),\ \forall \lambda\in (0,\delta_0).
\end{equation}
In fact, we shall prove that $\delta_0$ can be arbitrarily large in (\ref{F0G0-global-1}) which ends the proof of (\ref{F0G0-global}). By contradiction, let us assume that $
F_0(\lambda_0)<G_0(\lambda_0)\,
$ for some $\lambda_0>0$; clearly, $\lambda_0>\delta_0.$
Due to (\ref{F0G0-global-1}), we may set
\[
\lambda_s=\sup \{ \lambda<\lambda_0;\,\, F_0(\lambda)= G_0(\lambda) \}\,.
\]

\n Then, $\lambda_s<\lambda_0$ and for any $\lambda \in[\lambda_s,\lambda_0]$,
 one has $
F_0(\lambda)\leq G_0(\lambda).
$
For $\lambda>0$, we define the function
$\varphi_{\lambda}:(0,\infty)\rightarrow \R$  by
\[
\varphi_{\lambda}(t)=t^{\frac{n-4}{n}}-C(n-2)^2(n-4)^2\lambda^2t\,.
\]
We notice that  $\varphi_{\lambda}$ is non-decreasing in
$(0,t_{\lambda}]$, where
\[
t_{\lambda}=\frac{\lambda^{-\frac{n}{2}}}{(Cn(n-4)(n-2)^2)^{\frac{n}{4}}}\,.
\]
On one hand, a straightforward computation shows that for every $\lambda>0$, one has $$0<-\frac{ G_0'(\lambda)}{\lambda(n-2)(n-1)}=\left(\frac{K_0}{C}\right)^{\frac{n}{4}}\frac{ G''(\lambda)}{(n-2)(n-1)}<t_\lambda.$$
On the other hand, relation (\ref{masodrendu}) and the assumption $C\leq \frac{n+2}{n-2}K_0$ imply that  for every $\lambda>0$, 
$$0<-\frac{ F_0'(\lambda)}{\lambda(n-2)(n-1)}=\frac{ F''(\lambda)}{(n-2)(n-1)}\leq \frac{ G''(\lambda)}{(n-2)(n-1)}\leq t_\lambda.$$

We claim that 
\begin{equation}\label{F-G-monoton}
{ F_0'(\lambda)}\geq{ G_0'(\lambda)}\,,\ \forall \lambda \in [\lambda_s,\lambda_0].
\end{equation}
Since $
F_0(\lambda)\leq G_0(\lambda)
$ for every $\lambda \in[\lambda_s,\lambda_0]$,
by relations  (\ref{633}) and (\ref{6}) we have that 
\begin{eqnarray*}
	\varphi_{\lambda}\left(-\frac{ F_0'(\lambda)}{\lambda(n-2)(n-1)}\right)&=&
	\left( -\frac{ F_0'(\lambda)}{\lambda(n-2)(n-1)} \right)^{\frac{n-4}{n}}
	+ C(n-4)^2\frac{n-2}{n-1}\lambda F_0'(\lambda)\\
	&\leq &4C(n-4)^2F_0(\lambda)\\
	&\leq& 4C(n-4)^2G_0(\lambda)\\
	&=&\left( -\frac{ G_0'(\lambda)}{\lambda(n-2)(n-1)} \right)^{\frac{n-4}{n}}
	+ C(n-4)^2\frac{n-2}{n-1}\lambda G_0'(\lambda)\\&=&\varphi_{\lambda}\left(-\frac{ G_0'(\lambda)}{\lambda(n-2)(n-1)}\right),\ \ \ \ \forall \lambda \in[\lambda_s,\lambda_0].
\end{eqnarray*}
  By the monotonicity of $\varphi_\lambda$ on $(0,t_\lambda]$, relation (\ref{F-G-monoton}) follows at once. In particular, the function $F_0-G_0$ is non-decreasing on the interval  
$[\lambda_s,\lambda_0]$. Consequently,  we have
\[
0=F_0(\lambda_s)-G_0(\lambda_s)\leq F_0(\lambda_0)-G_0(\lambda_0)<0,
\]
 a contradiction, which shows the validity of  (\ref{F0G0-global}).

%
%
%
%
%
%
%
%
%
%
%
%
%
%

{\bf Step 5.} {\it Global volume non-collapsing property concluded.} 
 Inequality  (\ref{F0G0-global}) can be rewritten into 
\begin{equation}\label{9}
\int
_0^{\infty}\left({\rm vol}_g[B(x_0,t)]-b\omega_n t^n\right)\frac{((n-1)\lambda+t^2)t}{\left(
\lambda+t^2 \right)^{n}}{\rm d}t\geq 0\,,\ \forall \lambda>0,
\end{equation}
\n where
\[
b=(C^{-1}K_0)^{\frac{n}{4}}\,.
\]
The Bishop-Gromov comparison theorem  implies that the function
$t\mapsto \frac{{\rm vol}_g[B(x_0,t)]}{\omega_n t^n}$ is non-increasing on $(0,\infty)$; thus, the asymptotic volume growth 
\[
\limsup \limits_{t\rightarrow
\infty}\frac{{\rm vol}_g[B(x_0,t)]}{\omega_n t^n}=b_0\,
\]
is finite (and independent of the base point $x_0$).

  We shall prove that $b_0\geq b$. By contradiction, let us suppose that
$b_0=b-\varepsilon_0$ for some $\varepsilon_0>0$. Thus, there
exists a number $N_0>0$ such that
\begin{equation}\label{10}
\frac{{\rm vol}_g[B(x_0,t)]}{\omega_n t^n}\leq b-\frac{\varepsilon_0}{2}\,,
\,\,\forall t\geq N_0\,.
\end{equation}
For simplicity of notation, let 
$$f(\lambda,t)=\frac{((n-1)\lambda+t^2)t}{\left(
	\lambda+t^2 \right)^{n}},\ \ \lambda,t>0.$$
 Substituting (\ref{10}) into (\ref{9}) and by using the Bishop-Gromov comparison theorem, we obtain for every
$\lambda>0$ that
\begin{eqnarray*}
0&\leq& \int
_0^{\infty}\left({\rm vol}_g[B(x_0,t)]-b\omega_n t^n\right)f(\lambda,t){\rm d}t \\
&\leq & \int
_0^{N_0}{\rm vol}_g[B(x_0,t)]f(\lambda,t){\rm d}t+(b-\frac{\varepsilon_0}{2})\omega_n\int_{N_0}^\infty t^nf(\lambda,t){\rm d}t-b\omega_n\int_{0}^\infty t^nf(\lambda,t){\rm d}t\\
&\leq &\omega_n\int_{0}^{N_0} t^nf(\lambda,t){\rm d}t-b\omega_n\int_{0}^{N_0} t^nf(\lambda,t){\rm d}t-\frac{\varepsilon_0}{2}\omega_n\int_{N_0}^\infty t^nf(\lambda,t){\rm d}t\\&=&\omega_n(1-b+\frac{\varepsilon_0}{2})\int_{0}^{N_0} t^nf(\lambda,t){\rm d}-\frac{\varepsilon_0}{2}\omega_n\int_{0}^\infty t^nf(\lambda,t){\rm d}t.
\end{eqnarray*}
Note that for every $\lambda>0$, one has
\begin{eqnarray*}
I_1(\lambda)&=&\int_{0}^\infty t^nf(\lambda,t){\rm d}t=\lambda^\frac{4-n}{2}\int_{0}^\infty s^nf(1,s){\rm d}s\\
&=&\frac{2^{1-n}\pi^\frac{1}{2}(n^2-4n+6)\Gamma(\frac{n}{2}+1)}{(n-2)(n-4)\Gamma(\frac{n+1}{2})}\lambda^\frac{4-n}{2},
\end{eqnarray*}
and 
\begin{eqnarray*}
	I_2(\lambda)&=&\int_{0}^{N_0} t^nf(\lambda,t){\rm d}t=\int_{0}^{N_0} t^{n+1}\frac{(n-1)\lambda+t^2}{\left(
		\lambda+t^2 \right)^{n}}{\rm d}t\\
	&\leq&(n-1)N_0^{n+1}\lambda^{-n+1}+N_0^{n+3}\lambda^{-n}.
\end{eqnarray*}
Consequently, the above estimates show that for every $\lambda>0$,
$$M_0 \lambda^\frac{4-n}{2}\leq M_1\lambda^{-n+1}+M_2\lambda^{-n},$$
where $M_0,M_1, M_2>0$ are independent on $\lambda>0.$ It is clear that the latter inequality is not valid for large values of $\lambda>0$, i.e., we arrived to a contradiction. Accordingly, 
for every $r>0$, 
$$\frac{{\rm vol}_g[B(x_0,r)]}{\omega_n r^n}\geq \limsup \limits_{t\rightarrow
	\infty}\frac{{\rm vol}_g[B(x_0,t)]}{\omega_n t^n}=b_0\geq b=(C^{-1}K_0)^{\frac{n}{4}}.$$
Since the asymptotic volume
growth of $(M,g)$ is independent of the point $x_0$, we obtain the desired property, which completes the proof of  Theorem \ref{main}. \hfill $\square$

\begin{remark}\rm \label{remark-laplace}
	Note that relation (\ref{mainE2}) is equivalent to the  distance Laplacian growth condition. Indeed, a simple computation in Step 2 led us to relation (\ref{mainE2}) through the distance Laplacian growth condition. Conversely, if $\lambda\to 0$ in (\ref{mainE2}),  we obtain precisely that $\rho \Delta_g \rho \geq n-5.$ 
\end{remark}

{\it Proof of Theorem \ref{theorem-quantitative}.} 
(i) Due to Anderson \cite{Anderson} and Li \cite{Peter_Li},  if
 vol$_g[B(x,r)]\geq k_0 \omega_n
r^n$ for every $r>0$, then $(M,g)$ has finite fundamental
group $\pi_1(M)$ and its order is bounded above by ${k_0}^{-1}.$
By Theorem \ref{main} (ii) the property follows directly. In particular, if $C<2^\frac{4}{n}K_0$, then the order of $\pi_1(M)$ is strictly less than 2, thus $M$ is simply connected. 

(ii) First of all, due to Munn \cite[Table 5]{Munn-JGA} and a direct computation, for every $n\geq 5$ one has
$$\alpha_{MP}(1,n)^{-\frac{4}{n}}= 2^\frac{4}{n}<\frac{n+2}{n-2}.$$
Thus, since $\alpha_{MP}(\cdot,n)$ is increasing,  the values $\alpha_{MP}(k,n)^{-\frac{4}{n}}K_0$  are within the range where Theorem \ref{main} (ii) applies, $k\in \{1,...,n\}$. 

Now, let us assume that $ C <
\alpha_{MP}(k_0,n)^{-\frac{4}{n}}K_0$ for some $k_0\in
\{1,...,n\}$. By Theorem \ref{main} (ii) we have
the following estimate for the asymptotic volume
growth of $(M,g)$: 
$$\lim_{t\to \infty}\frac{{\rm vol}_g[B(x,t)]}{\omega_n t^n}\geq\left(\frac{K_0}{ C}\right)^\frac{n}{4}>\alpha_{MP}(k_0,n)\geq...\geq \alpha_{MP}(1,n).$$
Therefore, due to Munn \cite[Theorem 1.2]{Munn-JGA}, one has that
$\pi_1(M)=...=\pi_{k_0}(M)=0.$

(iii) If $C < \alpha_{MP}(n,n)^{-\frac{4}{n}}K_0$, then $\pi_1(M)=...=\pi_{n}(M)=0.$ Standard topological argument implies -based on Hurewicz's isomorphism theorem,- that $M$ is contractible.

(iv) If $ C = K_0$ then by Theorem
\ref{main} (ii) and the Bishop-Gromov volume
comparison theorem follows that vol$_g[B(x,r)]=\omega_nr^n$
for every $x\in M$ and $r>0$. Now, the equality in Bishop-Gromov
theorem implies that $(M,g)$ is isometric to the Euclidean space
$\mathbb R^n$. The converse is
trivial. \hfill $\square$

\section{Final remarks}\label{section-final}

\n We conclude the paper with some remarks and further questions: 

(a)   If $(M,g)$ is a complete $n-$dimensional Riemannian manifold and $x_0\in M$ is arbitrarily fixed, we notice that  $$\rho\Delta_g \rho=n-1+\rho\frac{J'(u,\rho)}{J(u,\rho)}\ \ {\rm a.e.\ on}\ M,$$
where $\rho(x)=\rho(x,x_0)$, $x=\exp_{x_0}(\rho(x)u)$ for some $u\in T_{x_0}M$ with $|u|=1$, and $J$ is the density of the volume form in normal coordinates, see Gallot,  Hulin and Lafontaine \cite[Proposition 4.16]{GHL}. On one hand, if the Ricci curvature on $(M,g)$ is nonnegative, one has $J'(u,\rho)\leq 0$. On the other hand,  the distance Laplacian growth condition 
$\rho \Delta_g \rho \geq n-5$ is equivalent to $$\frac{J'(u,\rho)}{J(u,\rho)}\geq -\frac{4}{\rho},$$ which is a curvature restriction on the manifold $(M,g)$. We are wondering if the latter condition can be removed from our results, which plays a crucial role in our arguments; see also Remark \ref{remark-laplace}. Examples of Riemannian manifolds  verifying the distance Laplacian growth condition (that are isometrically immersed into $\mathbb R^N$ with $N$ large enough) can be found in Carron \cite{Carron}.  

(b)  The requirement $C\leq \frac{n+2}{n-2}K_0$ is needed to explore the monotonicity of the function $\varphi_\lambda$ on $(0,t_\lambda]$, see Step 4 in the proof of Theorem \ref{main}. Although this condition is widely enough to obtain quantitative results, cf. Theorem \ref{theorem-quantitative}, we still believe that it can be somehow removed. 

(c) Let $(M,g)$ be an $n$-dimensional complete open Riemannian manifold 
with nonnegative Ricci curvature and fix $k\in \mathbb N$ such that $n>2k$. Let us consider for some $C>0$ the $k$-th order Sobolev inequality  
$$\left(\int_M |u|^{\frac{2n}{n-2k}}{\rm d}v_g\right)^\frac{n-2k}{n}\leq C \int_M (\Delta_g^{k/2} u)^2{\rm d}v_g,\ \forall u\in C^{\infty}_0(M),\eqno{({\bf SI})_C^k}$$
where 
\begin{equation*}
\Delta_g^{k/2} u= \left\{
\begin{array}{lll}
\Delta_g^{k/2}u & \mbox{if} & k \ {\rm is\ even}, \\
|\nabla_g(\Delta_g^{(k-1)/2}u)| & \mbox{if} &
k \ {\rm is\ odd}.
\end{array}%
\right.
\end{equation*}
Clearly,  $({\bf SI})_C^1=({\bf FSI})_C$ and  $({\bf SI})_C^2=({\bf SSI})_C$. It would be interesting to establish  $k$-th order counterparts of Theorems \ref{main}\&\ref{theorem-quantitative} with $k\geq 3$, noticing that the optimal Euclidean  $k$-th order Sobolev inequalities are well known with the optimal constant 
$$\Lambda_k=\left[\pi^k n(n-2k)\Pi_{i=1}^{k-1}(n^2-4i^2)\right]^{-1}\left(
\frac{\Gamma(n)}{\Gamma(\frac{n}{2})} \right)^{{2k}/{n}},$$
and the unique class of extremal functions (up to translations and multiplications) $$u_\lambda(x)=(\lambda+|x^2|)^\frac{2k-n}{2},\ x\in \mathbb R^n,$$
 see  Cotsiolis and Tavoularis \cite{CT}, Liu \cite{Liu-Calculus}. Once we use $w_\lambda=(\lambda+\rho^2)^\frac{2k-n}{2}$  as a test-function in 
 $({\bf SI})_C^k$, after a multiple application of the chain rule  we have to estimate in a sharp way the terms appearing in $\Delta_g^{k/2} w_\lambda$, similar to the eikonal equation $|\nabla_g \rho|=1$ and the distance Laplacian comparison  $\rho\Delta_g \rho \leq n-1$, respectively. In the second-order case this fact is highlighted in relation (\ref{mainE1}).  Furthermore,  higher-order counterparts of the distance Laplacian growth condition $\rho \Delta_g \rho \geq n-5$ should be found, (see relation (\ref{mainE2}) for the second order case), assuming this condition cannot be removed, see (a). 
%


%
%
%
%

%

\end{document}